%% file: nf-gan.tex
\documentclass[preprint,12pt]{elsarticle}

\usepackage{graphicx}
\usepackage{amssymb}
\usepackage{float}
\usepackage{lineno,hyperref}
\usepackage{amsmath,amsfonts,amsthm,amssymb,graphicx,subfigure,mathtools,tikz,hyperref,bm,blindtext}
\usepackage{diagbox}
\usepackage{geometry}
\usepackage[acronym]{glossaries}
\usepackage{multirow}
\usepackage{algorithm}
\usepackage{algpseudocode}
\input{math_commands}

\usepackage{lineno}

\journal{Journal Name}

\begin{document}

\begin{frontmatter}


\title{Variational inference in neural functional prior using normalizing flows: Application to differential equation and operator learning problems}



\author[hust]{Xuhui Meng\corref{cor}}
\cortext[cor]{Corresponding author}
\ead[cor]{xuhui\_meng@hust.edu.cn}

\address[hust]{Institute of Interdisciplinary Research for Mathematics and Applied Science, School of Mathematics and Statistics, Huazhong University of Science and Technology, Wuhan, 430074, China}

\begin{abstract}
Physics-informed deep learning have recently emerged as an effective tool for leveraging both  observational data and available physical laws. Physics-informed neural networks (PINNs) and deep operator networks (DeepONets) are two such models. The former encodes the physical laws via the automatic differentiation, while the latter learns the hidden physics from data. Generally, the noisy and limited observational data as well as the overparameterization in neural networks (NNs) result in uncertainty in predictions from deep learning models.  In \cite{meng2022learning}, a  Bayesian framework based on the {\emph{Generative Adversarial Networks}} (GAN) has been proposed as a unified model to quantify uncertainties in predictions of PINNs as well as DeepONets. 
Specifically, the proposed approach in \cite{meng2022learning} has two stages: (1) prior learning, and (2) posterior estimation. At the first stage, the GANs are employed to learn a functional prior either from a prescribed function distribution, e.g., Gaussian process, or from historical data and available physics. At the second stage,  the Hamiltonian Monte Carlo (HMC) method is utilized  to estimate the posterior in the latent space of GANs. However, the vanilla HMC does not support the mini-batch training, which limits its applications in problems with big data. In the present work, we propose to use the normalizing flow (NF) models in the context of variational inference, which naturally enables the minibatch training, as the alternative to HMC for posterior estimation in the latent space of GANs. A series of numerical experiments, including a nonlinear differential equation problem and a 100-dimensional Darcy problem, are conducted to demonstrate that NF with full-/mini-batch training are able to achieve similar accuracy as the ``gold rule'' HMC. Moreover, the minibatch training of NF makes it a promising tool for quantifying uncertainty in solving high-dimensional partial differential equation (PDE) problems with big data.
\end{abstract}

\begin{keyword}
uncertainty quantification \sep physics-informed neural networks \sep deep operator networks  \sep generative adversarial networks \sep normalizing flows \sep  differential equations 


\end{keyword}

\end{frontmatter}


\section{Introduction}
\label{sec:intro}
Physics-informed deep learning, capable of leveraging both data and physics, has emerged as an effective tool for diverse applications in the area of scientific computing, e.g., data fusion \cite{meng2019composite,meng2021fast,meng2021multi}, solving forward and inverse differential equations \cite{raissi2019physics,sirignano2018dgm,han2018solving,yu2018deep,lu2021learning,li2020fourier}.  Generally, the data noise associated with the measurement error in real-world applications and the overparametrization of NNs result in uncertainties in model predictions \cite{abdar2021review,pickering2022discovering,linka2022bayesian}. Quantifying the uncertainty propagation in deep learning approaches is crucial for critical applications involving physical and biological systems \cite{pickering2022discovering,linka2022bayesian}.

In the context of uncertainty quantification (UQ) in deep learning, the Bayesian neural network (BNN) has been one of the most popular and successful models for decades \cite{neal2012bayesian}. Recently, this approach has also been applied to solving forward and inverse ordinary/partial differential equation (ODE/PDE) problems \cite{yang2021b}, which is referred to as Bayesian physics-informed neural network (B-PINN) in \cite{yang2021b}. It is shown that B-PINN is capable of quantifying both the \emph{aleatoric} uncertainty associated with the data noise as well as the \emph{epistemic} uncertainty associated with the overparameterization of NNs. In addition to the BNN, extensions of standard neural networks, such as deep ensembles \cite{lakshminarayanan2017simple,pearce2020uncertainty} and dropout \cite{zhang2019quantifying}, has also been proposed to quantify uncertainties in deep learning. However, it is challenging for these models to achieve similar accuracy in terms of quality of the predicted uncertainties when comparing to BNNs \cite{yao2019quality,psaros2023uncertainty}. Interested readers are directed to \cite{yao2019quality,psaros2023uncertainty} for comprehensive comparisons among different UQ methods.

Despite the success of BNN, there are still limitations to be addressed. For instance, (1) the prior in BNN is specified for the weights/biases (i.e., hyperparameters), but the effect of the prior distribution for the hyperparamters on the prior in the functional space (i.e., output of BNNs) remains unclear, and (2) the number of dimensions for the hyperparameters in the BNN is generally  very high, resulting in difficulties in posterior estimation. To address these issues, \cite{meng2022learning} proposed to directly learn priors in functional space from data using generative adversarial networks (GANs) or physics-informed GANs (PI-GANs) if we have physical laws or PDEs. More specifically, the method proposed in \cite{meng2022learning} has two stages, i.e., (1) train a GAN/PI-GAN given data and/or physics to learn the prior in functional space, and (2) estimate the posterior distribution in the latent space of GAN/PI-GAN using Hamiltonian Monte Carlo (HMC) \cite{neal2012bayesian}. Note that the latent space of GAN/PI-GAN is generally characterized as low dimensional, which alleviates the difficulties in posterior estimation in BNNs/B-PINNs. It is reported that the approach proposed in \cite{meng2022learning} is capable of learning flexible functional priors, e.g., both Gaussian and non-Gaussian process, and it can also predict reasonable uncertainty bounds in regression as well as PDE problems. In the present work, we refer to the learned functional prior using GAN/PI-GAN as neural functional prior (NFP) since the prior is represented by NNs once the GAN/PI-GAN is well-trained.

Although HMC is capable of providing accurate posterior estimation as reported in \cite{meng2022learning}, the vanilla HMC does not support mini-batch training which restricts its applications in problems with big data. In addition to HMC, the variational inference (VI) is also a widely used approach for posterior estimation. In particular, the mean-field variational inference (MFVI), in which the true posterior is approximated by a Gaussian distribution, is one of the most commonly employed models \cite{blundell2015weight}. The MFVI is a unified method supporting both the full- and mini-batch training since it is based on stochastic gradient descent optimization algorithms. It can naturally handle problems with big data. However, as reported in \cite{yao2019quality,psaros2023uncertainty}, it is challenging for MFVI to achieve similar accuracy as the ``gold rule'' HMC in terms of the quality of the predicted uncertainties due to the fact that true posterior can be highly non-Gaussian. It has also been reported that employments of richer posterior approximations do result  in better performance in VI \cite{rezende2015variational}. Recently, the NN-based generative models, i.e., normalizing flows (NF), have been proposed to define  a wide range of probability distributions, e.g., Gaussian and non-Gaussian. Rezende {\sl et al.} \cite{rezende2015variational} utilized the NF to approximate the true posterior in VI, which clearly improves the accuracy compared to MFVI. Furthermore, the training of NF is also based on stochastic optimization as in the MFVI, indicating that it also enables the minibatch training.  In this study, we propose to employ the NF rather than HMC to compute the posterior in the latent space of NFP, which is expected to be a unified approach to enable both full- and mini-batch training, and thus making it a promising tool for solving high-dimensional parametric partial differential equation (PDE) problems with big data.

We organize the rest of this paper as follows. In Sec. \ref{sec:meth}, we review the neural functional prior for learning the priors from data and/or physics in functional space,  and we also introduce how to use NF to estimate the posterior in NFP. In Sec. \ref{sec:results}, we present results for three numerical experiments including regression and differential equation problems. We then summarize the findings of this study in Sec. \ref{sec:summary}. Finally, a brief overview for DeepONets and the details for the employed NNs as well as training strategy are provided in  \ref{sec:deeponet} and \ref{sec:nndetails}, respectively.

\section{Methodology}
\label{sec:meth}

Consider a nonlinear ODE/PDE problem for a certain physical system as follows: 
\begin{equation}\label{eqn:Nequ}
\begin{aligned}
        \mathcal{F}_{\lambda, \zeta}[u(\vx; \zeta)]& = f (\vx; \zeta), \vx \in \Omega, ~ \zeta \in \mathcal{Z}, \\
        \mathcal{B}_{\lambda, \zeta}[u(\vx; \zeta)]&=b(\vx; \zeta), \vx \in \partial\Omega
\end{aligned}
\end{equation}
where $\vx$ is the $D_x$-dimensional temporal-spatial coordinate,  $\zeta$ is a random event in a probability space $\mathcal{Z}$, $\mathcal{F}_{\lambda, \zeta}$ is a general differential operator, $\mathcal{B}_{\lambda, \zeta}$ is the boundary operator, $\Omega$ is the computational domain, $\partial\Omega$ is the temporal-spatial boundary, $f$ is the forcing term, $b$ represents the boundary/initial condition, and $\lambda$ are the problem parameters.

Our particular interest in the present study is focused on the following two scenarios, i.e., Case (I) \emph{inverse or mixed} ODE/PDE problems, where 
$\mathcal{F}$ are deterministic and specified. That is, $\zeta$ is fixed. $u$ and $f$ are partially known, $\mathcal{B}$ is either known or unknown, and 
$\lambda$ are either partially known or unknown. The objective is to predict $u$, $f$, as well as $\lambda$ in the entire domain $\Omega$; and Case (II) \emph{operator learning} problems, in which $\mathcal{F}$ are  stochastic and unknown (The stochasticity of $\mathcal{F}_{\lambda, \zeta}$ arises from the stochasticity in $\lambda(\bm{x};\zeta)$ with $\zeta \in \mathcal{Z}$), $\mathcal{B}$ and $\lambda$ can be either known or unknown. We would like to employ DeepONets to learn the operator mapping $f$ to $u$ if we have paired data on $(f, u)$, or we can  learn the operator mapping $(f, b, \lambda)$ to $u$ if we have paired data on $(f/b/\lambda, u)$. For simplicity, we will use $\lambda$ to represent all the inputs for DeepONets. In addition,  we assume that we have two types of data in all problems, i.e., historical data $\overline{\mathcal{D}}$ for prior learning, and testing data $\mathcal{D}$ for posterior estimation, similar as in \cite{meng2022learning,psaros2023uncertainty,zou2022neuraluq}.
 


\subsection{Prior learning in function space}
\label{sec:fp}

\begin{figure}
\centering
\includegraphics[width=0.9\textwidth]{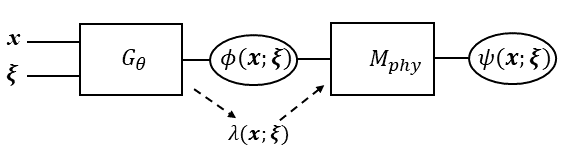}
\caption{\label{fig:nfp}
Schematic of the neural functional prior. $G_{\theta}$ is the generator in GANs parametized by $\theta$, which  takes $\bm{x}$ (spatial-temporal coordinates) and $\bm{\xi}$ (input noise in the latent space) as inputs; $\phi(\bm{x}, \bm{\xi})$ and $\psi(\bm{x}, \bm{\xi})$ can be the solution $u$ and the source $f$ in ~\Eqref{eqn:Nequ}, or vice-versa; $M_{phy}$ represents the available physics, which is either the differential operators in a certain ODE/PDE implemented by the automatic differentiation, or a pretrained DeepONet that encodes the hidden physics. Finally, $\tilde{\lambda}(\bm{x}, \bm{\xi})$ are either unknown or partially known problem parameters in inverse problems. Adapted from \cite{zou2022neuraluq,psaros2023uncertainty}.}
\end{figure}

In this section, we introduce how to learn the functional prior using GANs given historical data $\overline{\mathcal{D}}$ for solving problems of interest in \Eqref{eqn:Nequ}.

For Case (I), we employ the generator of GANs as the surrogate for $u$, which takes as inputs the coordinates $\bm{x}$ and the samples $\bm{\xi}$ from a prescribed distribution $P(\bm{\xi})$, e.g., uniform or Gaussian distribution. We refer to $\bm{\xi}$ as input noise in the present study. Also,  the standard multivariate Gaussian distribution is utilized for $P(\bm{\xi})$ in this work. As shown in Fig. \ref{fig:nfp}, we denote the generator as $G_{\theta}$, where $\theta$ represents the parameters in NNs, and the output of the generator is represented by $\tilde{u}(\bm{x}, \bm{\xi})$. For problems with multiple fields, we can either use multiple generators with one output each or a generator with multiple outputs as the corresponding surrogates. With specified $\mathcal{F}$, we can then obtain the surrogate for $f$, i.e., $\tilde{f}(\bm{x}, \bm{\xi}) = \mathcal{F}[\tilde{u}(\bm{x}, \bm{\xi})]$, using the automatic differentiation as in PINNs \cite{raissi2019physics}. We note that we can get the surrogate for $b$ in the similar way as obtaining $\tilde{f}$ if $\mathcal{B}$ is specified. For the inverse problems we consider here, we also need surrogates for the problem parameters, i.e., $\lambda$. In the current study, we employ another generator with the output $\tilde{\lambda}(\bm{x}; \bm{\xi})$ to approximate $\lambda$. 



We now discuss the historical data for training GANs. On the one hand, we have a certain number of samples  from a hidden distribution $P_r$, e.g., $\mathcal{\overline{D}} = \left[\mathcal{\overline{D}}_{u_i} \cup \mathcal{\overline{D}}_{f_i} \cup \mathcal{\overline{D}}_{b_i} \cup \mathcal{\overline{D}}_{{\lambda}_i}\right]^{\overline{N}}_{i=1}$, where $\overline{N}$ is the number of samples for $u/f/b/\lambda$. With the generators $\tilde{u}(\vx;\bm{\xi})$, $\tilde{f}(\vx;\bm{\xi})$, $\tilde{b}(\vx;\bm{\xi})$, and $\tilde{\lambda}(\vx;\bm{\xi})$, we can then train $G_{\theta}(\bm{x}, \boldsymbol{\xi})$ to generate ``fake'' samples that match samples from the hidden distribution $P_r$. In particular, we employ a certain numbers of discrete points to resolve each sample in $\overline{\mathcal{D}}$. Details will be presented in each test case of Sec. \ref{sec:results}. In addition, we denote the discriminator neural network in GANs by $D_{\bm{\eta}}$, which takes a real or fake sample as input, and outputs a real number. The loss function for optimizing the generator parameters $\bm{\theta}$ and the discriminator parameters $\bm{\eta}$ are as follows:
\begin{equation} \label{eqn:GDloss}
\begin{aligned}
L_G &= -\mathbb{E}_{\bm{\xi}\sim P(\bm{\xi})} [ D_{\bm{\eta}}(G_{\bm{\theta}}(\bm{x}, \boldsymbol{\xi}))], \\
L_D &= \mathbb{E}_{\bm{\xi}\sim P(\bm{\xi})} [ D_{\bm{\eta}}(G_{\bm{\theta}}(\bm{x}, \boldsymbol{\xi}))] - \mathbb{E}_{T\sim P_r} [ D_{\bm{\eta}}(T)] + \omega \mathbb{E}_{\hat{T} \sim P_i} (\Vert \nabla_{\hat{T}} D_{\bm{\eta}}(\hat{T}) \Vert_2 -1 )^2.
\end{aligned}
\end{equation}
where $P_i$ is the distribution induced by uniform sampling on interpolation lines between independent samples of $T$ and $G_{\bm{\theta}}(\bm{x}, \bm{\xi})$, and $\omega$ is the gradient penalty coefficient. Here we set $\omega=0.1$. During the training, we update $\bm{\eta}$ and $\bm{\theta}$ iteratively with the ratio of $5:1$ \cite{meng2022learning}. 


For Case (II), i.e., operator learning problems, we assume that we have paired data on $\lambda$ and $u$, we can then train a DeepONet to approximate the operator mapping $\lambda$ to $u$, which is denoted by $u = S[\lambda](\bm{x})$. A brief overview of DeepONets is presented in \ref{sec:deeponet}. Once the DeepONet is well-trained, it is capable of encoding the hidden physics represented by data. Afterwards, we employ the generator $\tilde{\lambda}(\bm{x}, \bm{\xi})$ in GANs as the surrogate for $\lambda$, and perform the training using the historical data on $\lambda$ or the same data for training the DeepONet, in a similar way as in Case (I).  Upon the completion of training GANs, we can obtain the functional priors for both $\lambda$ and $u$, which are denoted by $\tilde{\lambda}(\bm{x}, \bm{\xi})$ and $\tilde{u}(\bm{x}, \bm{\xi}) = S[\tilde{\lambda}(\bm{x}, \bm{\xi})] (\bm{x})$, respectively.

\subsection{Posterior estimation using normalizing flows}\label{sec:post}
Once the GANs are well-trained, we can then obtain $\tilde{u}(\bm{x}, \bm{\xi})$, and $\tilde{f}(\bm{x}, \bm{\xi})$/$\tilde{b}(\bm{x}, \bm{\xi})$/$\tilde{\lambda}_{\bm{\eta}}(\vx;\bm{\xi})$ represented by NNs to serve as the functional priors for new tasks in the future. We refer to the functional priors here as neural functional priors because we employ NNs to represent them. With the neural functional priors, the objective is to infer the posterior distribution in the latent space of GANs, i.e., $\bm{\xi}$, given data on a new task, i.e., $\mathcal{D}$. Specifically, the posterior for $\bm{\xi}$ can be expressed as follows based on the Bayes' rule:
\begin{equation}\label{eqn:Bayes}
\begin{aligned}
P(\bm{\xi}| \mathcal{D}) = \frac{ P(\mathcal{D}|\bm{\xi})P(\bm{\xi})}{P(\mathcal{D})},
\end{aligned}
\end{equation}
where $P(\bm{\xi})$ is the prior distribution for $\bm{\xi}$, which is Gaussian, i.e, $P(\bm{\xi}) = (2\pi)^{-d_{\bm{\xi}}/2}\exp(-\Vert \bm{\xi} \Vert^2/2)$ ($d_{\bm{\xi}}$ is the dimensionality of $\bm{\xi}$), here. In addition, $P(\mathcal{D}|\bm{\xi})$ is the likelihood for the testing data on the new task, i.e., $\mathcal{D} = \mathcal{D}_u \cup \mathcal{D}_f \cup \mathcal{D}_b \cup \mathcal{D}_{\lambda} $.  We can then write the likelihood as follows if we assume the measurement errors are Gaussian:
\begin{equation}
\label{eqn:ganlikelihood}
\begin{aligned}
    P(\mathcal{D}|\boldsymbol{\bm{\xi}}) &= P(\mathcal{D}_u|\bm{\xi}) P(\mathcal{D}_f|\bm{\xi}) P(\mathcal{D}_b|\bm{\xi}) P(\mathcal{D}_{\lambda}|\bm{\xi}), \\
     P(\mathcal{D}_{\alpha}|\bm{\xi}) &= \prod_{i=1}^{N_{\alpha}} \frac{1}{\sqrt{2\pi{\sigma_{\alpha}^{(i)}}^2}}\exp \left(-\frac{(\tilde{{\alpha}}(\boldsymbol{x}_{{\alpha}}^{(i)}; \bm{\xi}) - {\alpha}^{(i)})^2}{2{\sigma_{\alpha}^{(i)}}^2}\right), \alpha = u, f, b, \lambda, \\
\end{aligned}
\end{equation}
where $N_{\alpha}$ and $\bm{x}_{\alpha}$ are the number of discrete points and the corresponding coordinate for the field of $u/f/b/\lambda$, respectively; and $\sigma_{\alpha}$ is the standard deviation for the noise of the corresponding filed. Finally, the marginal likelihood $P(\mathcal{D})$ in \Eqref{eqn:Bayes} is in general analytically intractable, we will thus compute $P(\bm{\xi}|\mathcal{D})$ numerically using variational inference with normalizing flows. 



\begin{figure}
\centering
\includegraphics[width=0.8\textwidth]{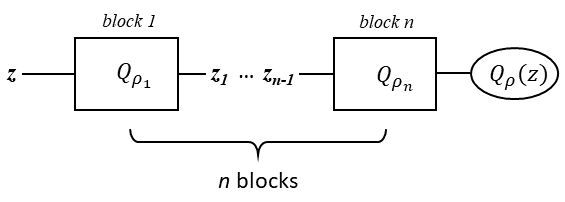}
\caption{\label{fig:nf}
Schematic of the normalizing flows. The NF consists of $n (n \ge 1)$ blocks, and each block is a NN. $\bm{z}$ denotes the input of NF, which is a multivariate standard normal distribution. In addition, $\bm{z}$ has the same dimensionality as $\bm{\xi}$. $\bm{z}_n$ is the output of $n_{th}$ block; $\rho_n$ represents the hyperparameters in the NN of $n_{th}$ block, and $G_{\rho}(\bm{z})$ is the final output of NF, where $\rho = (\rho_1, ..., \rho_n)$ is the collection of hyperparameters in all blocks. 
}
\end{figure}

In the variational inference,  the posterior density function for the unknown parameter  $\bm{\xi} = (\xi_1, \xi_2...\xi_{d_{\bm{\xi}}})$, i.e., $P(\bm{\xi} | \mathcal{D})$, is approximated by another density function, which is from a normalizing flow model parameterized by $\rho$, as illustrated in Fig. \ref{fig:nf}. We denote the samples and the corresponding density function from NF as $z_n = G_{\rho}(\bm{z})$ and $Q_{\rho}(\bm{z}_n) = Q_{\rho}(G_{\rho}(\bm{z}))$, respectively. We can then  tune $\rho$ to minimize
\begin{equation}\label{eqn:vi_nf}
    \begin{aligned}
    D_{KL}(Q_{\rho}(\bm{z}_n)||P(\bm{\xi}|\mathcal{D})) & \simeq \mathbb{E}_{z_n \thicksim Q_{\rho}(z_n)} [\ln Q_{\rho}(z_n) - \ln P(\bm{z_n}) - \ln P(\mathcal{D} | \bm{z}_n)] \\
    & = \frac{1}{N_z} \sum^{N_z}_{j=1} [\ln Q_{\rho}({\bm{z}_{n,j}}) - \ln P(\bm{z}_{n,j}) - \ln P(\mathcal{D} | \bm{z}_{n,j})],
    \end{aligned}
\end{equation} 
where $D_{KL}$ denotes the Kullback-Leibler (KL) divergence, and $N_z$ is the number of posterior samples for $\bm{\xi}$ used to compute the loss at each training step. Similar as in \Eqref{eqn:ganlikelihood}, $\mathcal{D}$ represents all available measurements on $u/f/b/\lambda$. Further, the loss function can be written as follows if the minibatch training is used:
\begin{equation}\label{eqn:vi_nf_mini}
    \begin{aligned}
    D_{KL}(Q_{\rho}(\bm{z}_n)||P(\bm{\xi}|\mathcal{D})) \simeq & \frac{1}{N_z} \sum^{N_z}_{j=1} [\ln Q_{\rho}(\bm{z}_{n,j}) - \ln P(\bm{z}_{n,j}) ~ -   
    \frac{N_{\alpha}}{M_{\alpha}}\sum^{M_{\alpha}}_{i=1}\ln P(\mathcal{D}_{{\alpha}_i} | \bm{z}_{n,j}) ], \alpha = u, f, b, \lambda, 
    \end{aligned}
\end{equation} 
where $N_{\alpha}$ and $M_{\alpha}$ are the numbers of all measurements and the minibatch size for the field $u/f/b/\lambda$, respectively. The detailed algorithm for VI with normalizing flows is presented in Algorithm \ref{alg:vi}. Note that Various variants of NF \cite{dinh2016density,papamakarios2017masked,kingma2016improved,kingma2018glow} have been developed recently. Here we employ the inverse autoregressive flow model proposed in \cite{kingma2016improved}  due to its efficiency in generating posterior samples.

\begin{algorithm}[H]
\caption{Variational inference with normalizing flows}
\label{alg:vi}
\begin{algorithmic}
\Require Pretrained  neural functional prior.
\Require An initialization for $\rho$.
\For{$k=1,2...N$}
\State Sample $\{\bm{z}^{(j)}\}_{j=1}^{N_z}$ independently from $\mathcal{N}(\bm{0}, \bm{I}_{d_{\xi}})$, \;
\State Compute the loss, i.e., $L(\rho)$, based on \Eqref{eqn:vi_nf} or \Eqref{eqn:vi_nf_mini}, \;
\State Update $\rho$ with gradient $\nabla_{\rho}L(\rho)$ using Adam optimizer. \;
\EndFor
\State Sample $\{\boldsymbol{z}^{(j)}\}_{j=1}^{M}$ independently from $\mathcal{N}(\bm{0}, \bm{I}_{d_{\bm{\xi}}})$, \;
\State Obtain the posterior samples using $\{\boldsymbol{z_n}^{(j)}\}_{j=1}^{M} = G_{\rho}(z)$, \;
\State Calculate $\{ \tilde{u}(\boldsymbol{x},\boldsymbol{z_n}^{(j)})\}_{j=1}^M$ as posterior samples for $u(\boldsymbol{x})$, similarly for other terms. \;

\end{algorithmic}
\end{algorithm}

With the posterior samples $\{\tilde{u}(\boldsymbol{x}, {\bm{z}}_n^{(j)} ) \}_{j=1}^M$ for $u$ in Algorithm \ref{alg:vi},  we can obtain the quantities of interest, i.e.,  the mean and standard deviation of $\{\tilde{u}(\boldsymbol{x}, {\bm{z}}_n^{(j)} ) \}_{j=1}^M$. The former represents the prediction of  $u(\boldsymbol{x})$ while the latter quantifies the uncertainty. In the present work, we set $M = 1,000$. Similarly for the other terms, e.g., $f/b$.
\section{Results and discussion}
\label{sec:results}
In this section, we first employ the normalizing flows for posterior estimation in the neural functional prior using the example of a one-dimensional function approximation case. We then test two differential equation problems in conjunction with PINN as well as DeepONet. In each case, we demonstrate that NF supports the minibatch training in posterior estimation. 


\subsection{Pedagogical example: Function Approximation}
\label{sec:func}

We first consider to employ the NFP as well as the NF to quantify uncertainties in a one-dimensional regression problem. The target function is expressed as follows:
\begin{equation}
 u = \sin^3(3x), x \in [-1, 1].
\end{equation}
We assume that we have 128 noisy measurements on this function which serves as the training data. Further, the training data are equidistantly distributed in $x \in [-0.8, -0.2] \cup [0.2, 0.8]$. The measurement error is assumed to be a Gaussian distribution with zero mean and 0.1 as the standard deviation.

For this specific case, we assume that the historical data are from the following Gaussian process, i.e.,
\begin{equation}
\begin{split}\label{eq:gp_prior}
    \bm{u} & \thicksim \mathcal{GP}(0, \mathcal{K}), ~ \mathcal{K} = \exp \left( - \frac{(x - x')^2}{2 l^2}\right),\\
    x,&~ x' \in [-1, 1], ~l = 0.2.
\end{split}
\end{equation}
We then randomly draw 10,000 samples from the aforementioned GP to train the GANs to obtain the prior in the functional space. In addition, 30 equidistant discrete points are employed to resolve each sample. The details for the architecture as well as the training of GANs are presented in \ref{sec:nndetails}, which will not be presented here. 

\begin{figure}
\centering
\subfigure[]{\label{fig:func_a}
\includegraphics[width=0.33\textwidth]{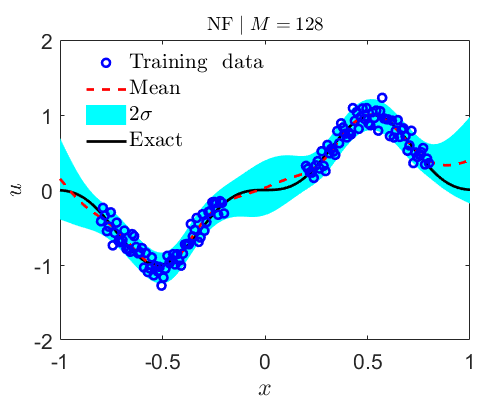}
\includegraphics[width=0.33\textwidth]{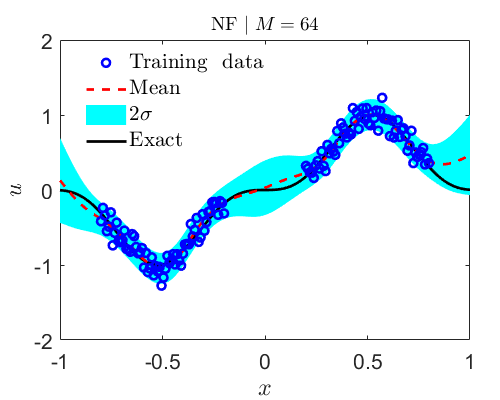}
\includegraphics[width=0.33\textwidth]{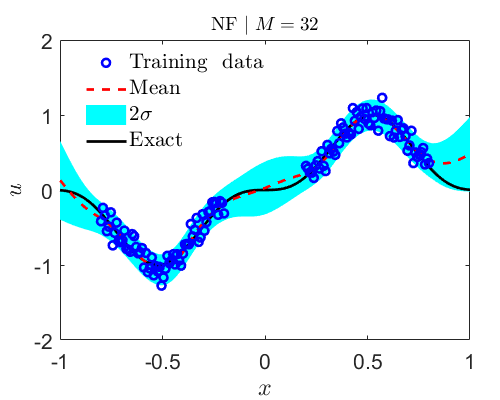}
}
\subfigure[]{\label{fig:func_b}
\includegraphics[width=0.33\textwidth]{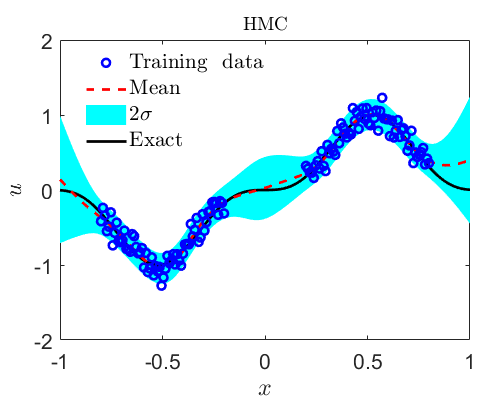}
}
\caption{\label{fig:func}
1D function approximation with 128 training data: Predictions from  (a) NF with full batch size (left most), batch size 64 (middle), and 32 (right most); and (b) HMC. $2 \sigma$: two standard deviations.
}
\end{figure}

With the learned neural functional prior, we now employ the NF to estimate the posterior distribution of $\bm{\xi}$ in the neural functional prior given testing data on an unseen task. We first conduct the fullbatch training in NF, i.e., $M = 128$, and
illustrate the results in Fig. \ref{fig:func_a}. As shown, (1) the predicted uncertainty increases at the regions that do not have training data, and (2) the computational errors between the predicted mean and the reference solution are bounded by the predicted uncertainty. To further test the minibatch training in NF, we then train the NF with two different batch sizes, i.e., $M = 32$ and 64. The results are also illustrated in Fig. \ref{fig:func_a}. It is observed that (1) both the predicted mean and uncertainty from NF with the two batch sizes are similar, and (2) the results in these two cases agree well with those in the first test case, i.e., NF with fullbatch training, suggesting the effectiveness of minibatch training in NF.

Finally, the results from the ``gold rule'', i.e., HMC, are utilized as the reference in Fig. \ref{fig:func_b}. Note that the minibatch training is not supported in vanilla HMC, we therefore only conduct the HMC with fullbatch training in the present study. As shown, both the predicted mean and uncertainty of NF with full- and mini-batch training are similar as those from HMC, suggesting that NF is able to achieve similar accuracy as HMC. The minibatch training in NF makes it a promising tool for quantifying uncertainty in problems with big data, which outperforms the HMC.

\subsection{PINNs: 1D nonlinear diffusion-reaction problem}
\label{sec:poisson}

We now consider to employ the neural functional prior and the NF to quantify uncertainties in an inverse differential equation problem. Specifically, the PI-GAN is utilized to learn the functional prior from historical data as well as the equation, and the NF is used to obtain the posterior samples in the latent space of PI-GAN given testing data. Specially, we consider the same case as in \cite{meng2022learning}, i.e., a nonlinear diffusion-reaction system, which is governed by: 
\begin{equation}\label{eq:reac}
\begin{split}
    D \partial^2_x u - k_r u^3 &= f, ~ x \in [-1, 1], \\
    u(-1) &= u(1) = 0,
\end{split}
\end{equation}
where $u$ is the solute concentration, $D = 0.01$ is the diffusion coefficient, $k_r$ is the chemical reaction rate, and $f$ is the source term. The exact solution for this system is assumed to be as follows:
\begin{align}\label{eq:reac_sol}
    u = (x^2 - 1) \sum^4_{i=1} \left[\omega_{2i-1} \sin(i \pi x) + \omega_{2i} \cos(i \pi x) \right],
\end{align}
where $\omega_{i}$ are uniformly sampled from $\mathcal{U}([0,1])$, $i = 1,2,...,8$. The chemical reaction rate is set as a nonlinear function of the solute concentration, i.e.,  $k_r = 0.4 \exp(-u)$, and $f$ can then be derived from  \Eqref{eq:reac}. 

\begin{figure}[H]
\centering
\subfigure[]{\label{fig:1dinversea}
\includegraphics[width=0.33\textwidth]{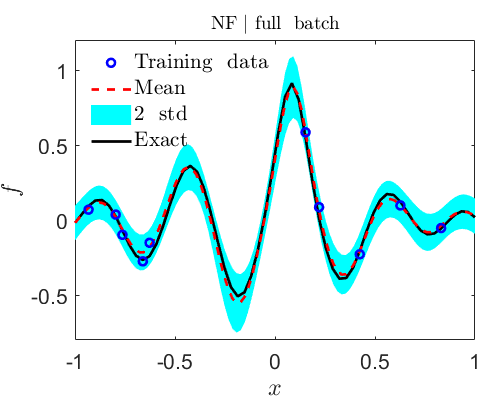}
\includegraphics[width=0.33\textwidth]{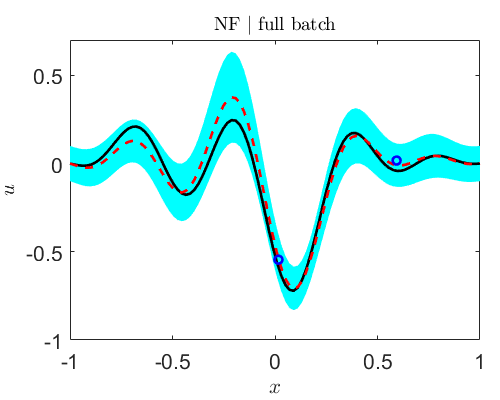}
\includegraphics[width=0.33\textwidth]{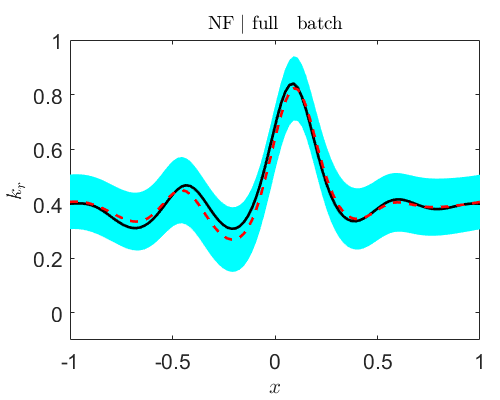}
}
\subfigure[]{\label{fig:1dinverseb}
\includegraphics[width=0.33\textwidth]{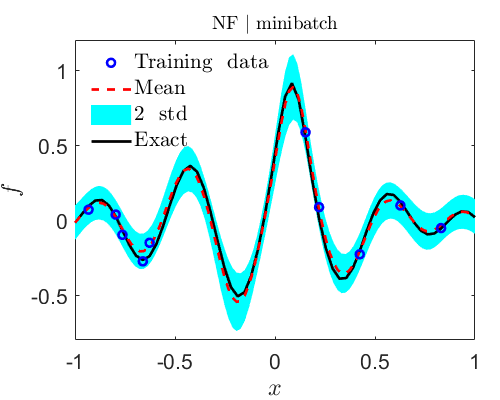}
\includegraphics[width=0.33\textwidth]{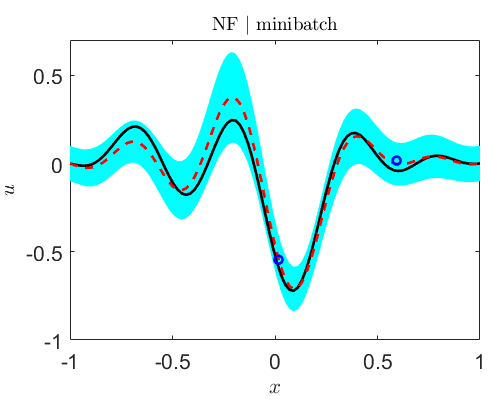}
\includegraphics[width=0.33\textwidth]{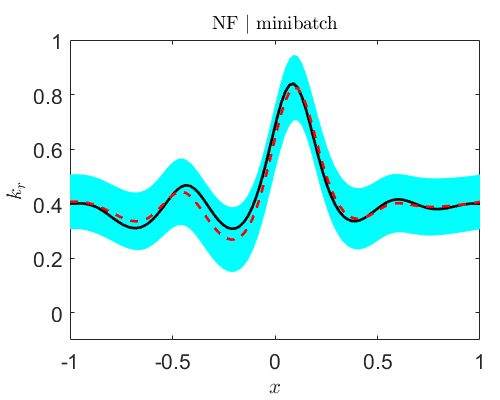}
}
\subfigure[]{\label{fig:1dinversec}
\includegraphics[width=0.33\textwidth]{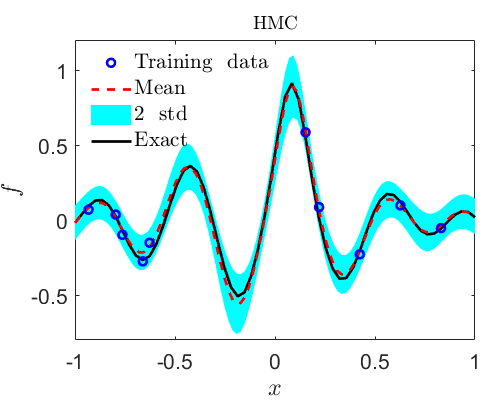}
\includegraphics[width=0.33\textwidth]{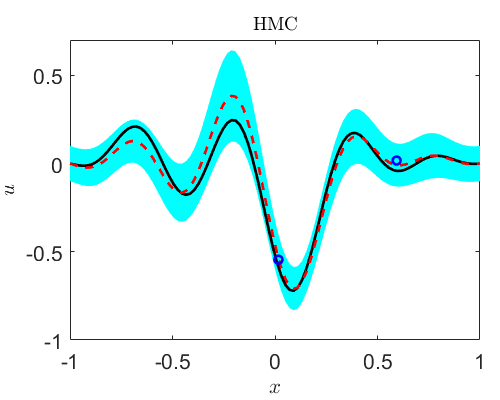}
\includegraphics[width=0.33\textwidth]{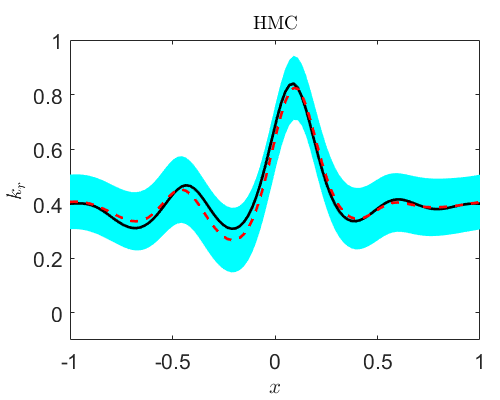}
}
\caption{\label{fig:1dinverse}
1D nonlinear diffusion-reaction problem: Predicted $f$, $u$ and $k_r$ from (a) NF with full batch training, (b) NF with minibatch training in which $M_u = 1$ and $M_f = 5$, and (c) HMC.
}
\end{figure}

Similar as in \cite{meng2022learning}, we assume that we have 10,000 pairs of $(k_r, f)$ as the historical data for learning the functional priors of $k_r$ and $f$.  Measurements from 40 equidistant sensors are used to resolve each $k_r/f$ sample. For the inverse problem we consider here, we assume that we have partial  noisy measurements on $u$ and $f$ for a new task,  the objective is to infer $k_r$ with uncertainties given data on $u$ and $f$. Specifically, we utilize 10 and 2 measurements for $f$ and $u$ as the training data at the posterior estimation stage, respectively. We note that the historical data as well as the setup, and the training data for the inverse problem are the same as in \cite{meng2022learning}.

We first train a PI-GAN to obtain the neural functional priors for $k_r$ and $f$ based on the historical data (Details for the architectures of PI-GAN can be found in \ref{sec:nndetails}). With the neural functional priors, we then employ the NF with full- and mini-batch training to compute the posterior given noisy measurements on $u/f$ in the new task. The results from full- and mini-batch training are illustrated in Figs. \ref{fig:1dinversea}-\ref{fig:1dinverseb}, respectively. As shown, (1) the predictions for $f$, $u$, and $k_r$ are quite similar, (2) the predicted $k_r$ agree well with the reference solution even we do not have any measurements on it, which is attributed to the informative prior learned from historical data, and  (3) the computational errors between the predicted means and the reference solutions for $f$, $u$ and $k_r$ are bounded by the predicted uncertainties. Finally, we present the results from the HMC to serve as the reference solutions in Fig. \ref{fig:1dinversec}. As observed, the results from NF with full-/mini-batch training are similar as those from HMC, which demonstrate that NF with full-/mini-batch training is able to achieve similar accuracy as HMC in posterior estimation.  

\subsection{DeepONet: 100-dimensional Darcy problem}
\label{sec:darcy}

We here employ the neural functional prior to quantify uncertainties in the predictions of DeepONet. Particularly, we utilize the test case in \cite{psaros2023uncertainty,zou2022neuraluq}, which is
a problem of two-dimensional steady flow through porous media. The governing equation for this problem is described by Darcy's law as follows:
\begin{linenomath*}
\begin{align}\label{eq:comp:don:pde}
	\nabla \cdot (\lambda(x, y) \nabla u(x, y)) =  f, ~ 0 \le x, y \le 1, 
\end{align}
\end{linenomath*}
where $x, y$ are the space coordinates,  $u(x, y)$ is the hydraulic head, $f$ is a constant, i.e., $f = 50$, and $\lambda(x, y)$ denotes the hydraulic conductivity field. The boundary conditions are expressed as follows:
\begin{linenomath*}
\begin{subequations}\label{eq:comp:don:bcs}
	\begin{align}
		u(x=0, y; \xi) = 1, ~ u(x=1, y; \xi) = 0, \\
		\partial_{\boldsymbol{n}} u(x, y=0; \xi) = \partial_{\boldsymbol{n}} u(x, y=1; \xi) = 0, \forall \xi \in \Xi,
	\end{align}
\end{subequations}
\end{linenomath*}
where $\boldsymbol{n}$ denotes the unit normal vector of the boundary. Generally, $\lambda(x, y)$ is determined by the pore structure. We then employ 
 a stochastic model for $\lambda(x, y; \zeta)$, to take in account of different porous media. In particular, $\lambda(x, y; \zeta) = \exp(\bar{\lambda}(x, y; \zeta))$, where $\bar{\lambda}(x, y; \zeta)$ is a truncated Karhunen-Lo\`{e}ve expansion of a GP with zero mean and kernel given by
\begin{linenomath*}
\begin{subequations}\label{eq:comp:don:gp}
	\begin{align}
		k_{\tilde{\lambda}}(x, y, x', y') = \exp(-\frac{(x - x')^2}{2l^2} - \frac{(y - y')^2}{2l^2}), \\
		0 \le x, y, x', y' \le 1, ~ l = 0.25. 
	\end{align}
\end{subequations}
\end{linenomath*}
In the following, we only keep the first 100 leading terms of the expansion and refer to the current problem as a 100-dimensional Darcy problem here.

\begin{figure}
\centering
\subfigure[]{
\includegraphics[width=0.3\textwidth]{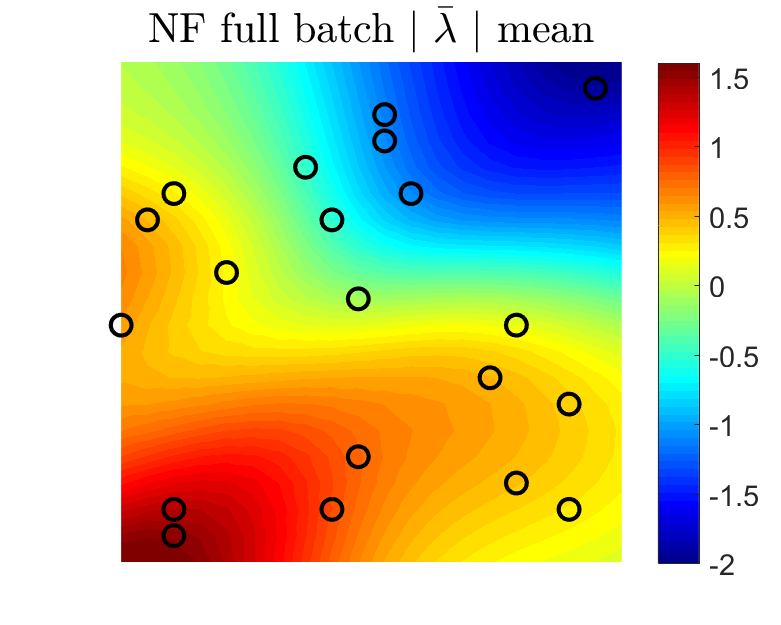}
\includegraphics[width=0.3\textwidth]{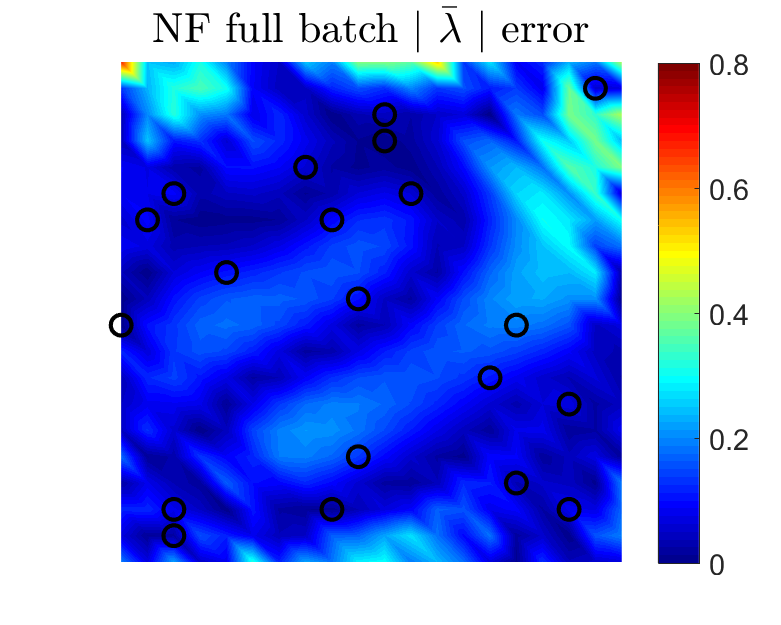}
\includegraphics[width=0.3\textwidth]{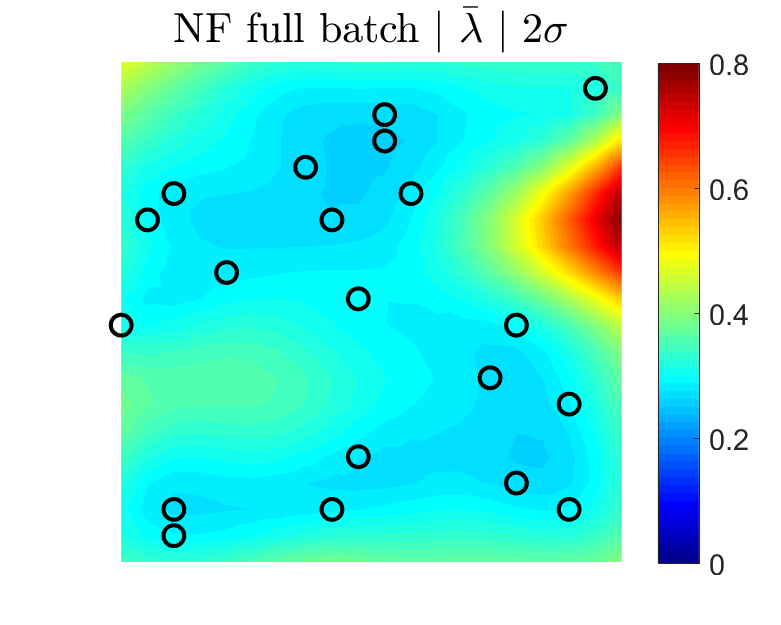}
}
\subfigure[]{
\includegraphics[width=0.3\textwidth]{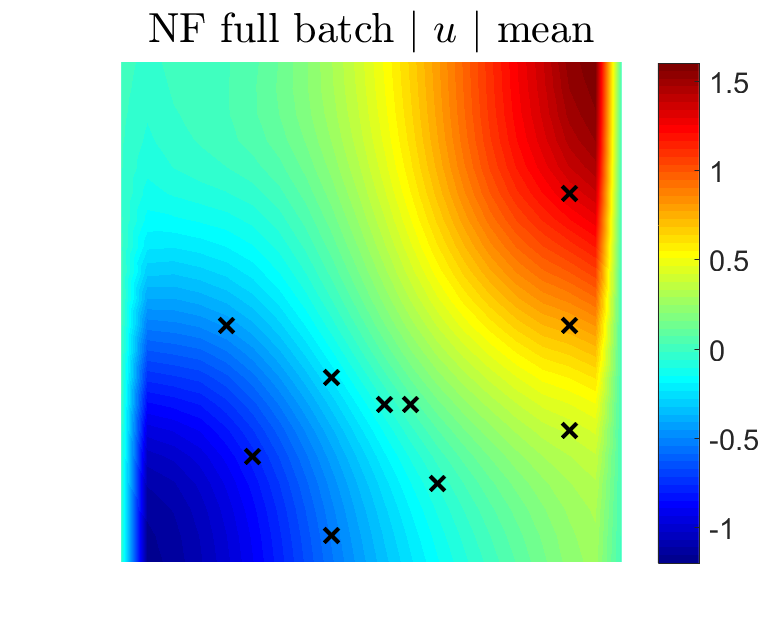}
\includegraphics[width=0.3\textwidth]{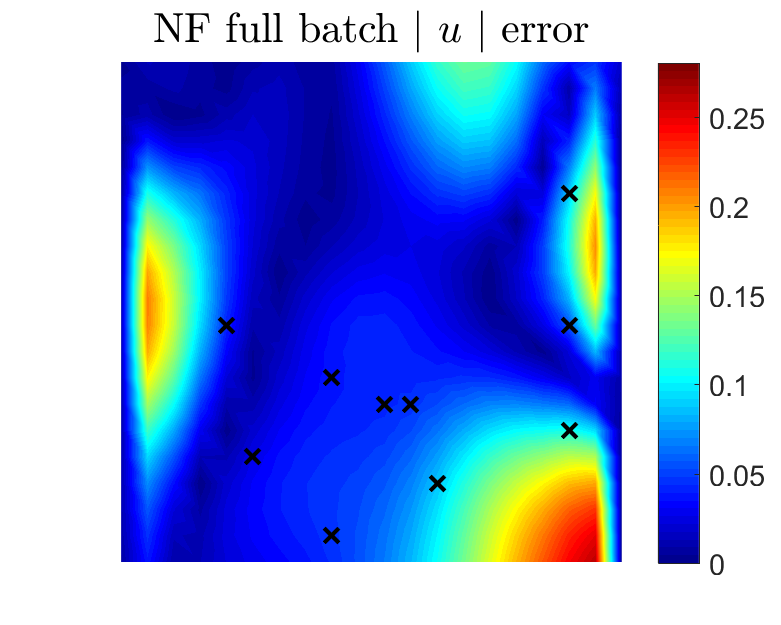}
\includegraphics[width=0.3\textwidth]{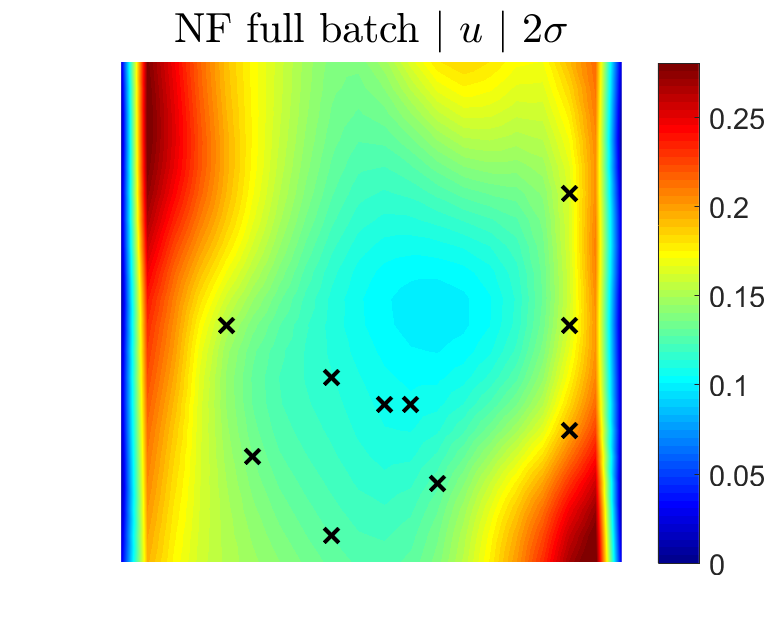}
}
\caption{\label{fig:darcy_nf_fb}
100D Darcy problem: Predictions from NF with full batch training. (a) $\bar{\lambda}$. (b) $u$. Circle: testing data for $\bar{\lambda}$, Cross: testing data for $u$.
}
\end{figure}

\begin{figure}
\centering
\subfigure[]{
\includegraphics[width=0.3\textwidth]{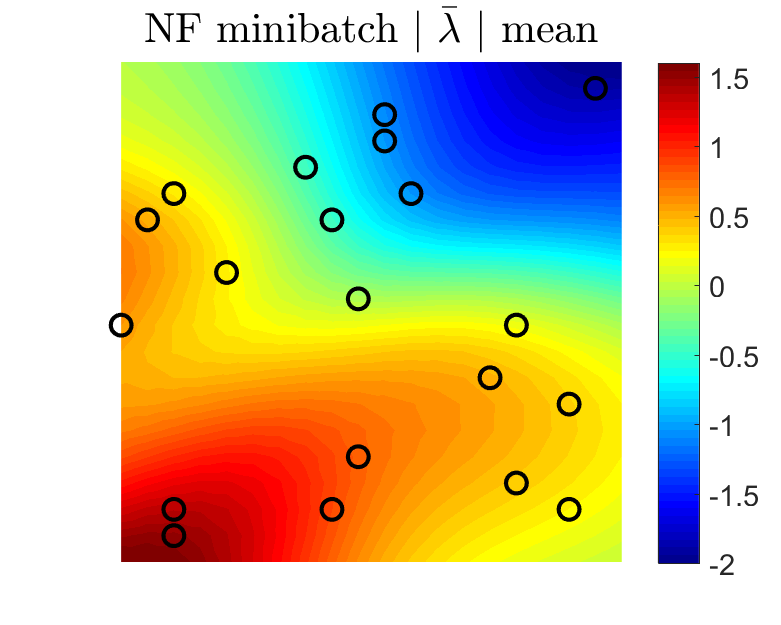}
\includegraphics[width=0.3\textwidth]{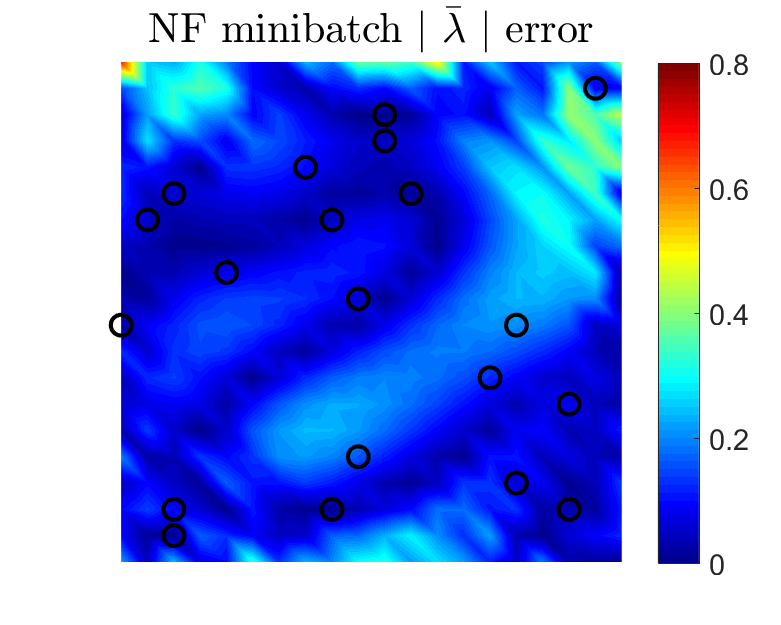}
\includegraphics[width=0.3\textwidth]{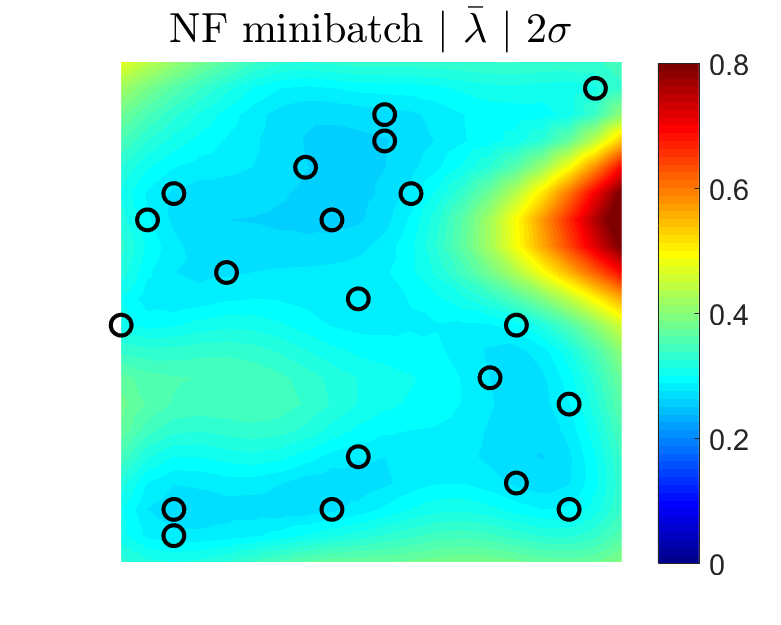}
}
\subfigure[]{
\includegraphics[width=0.3\textwidth]{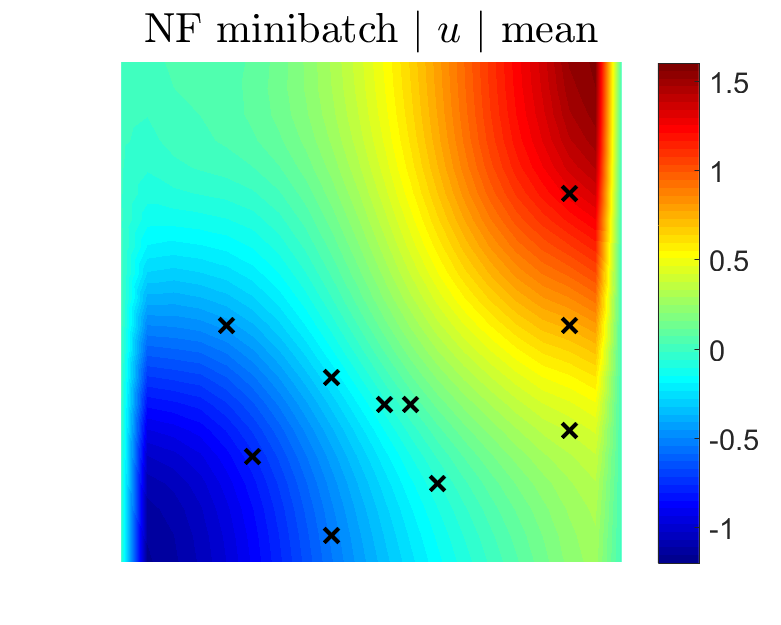}
\includegraphics[width=0.3\textwidth]{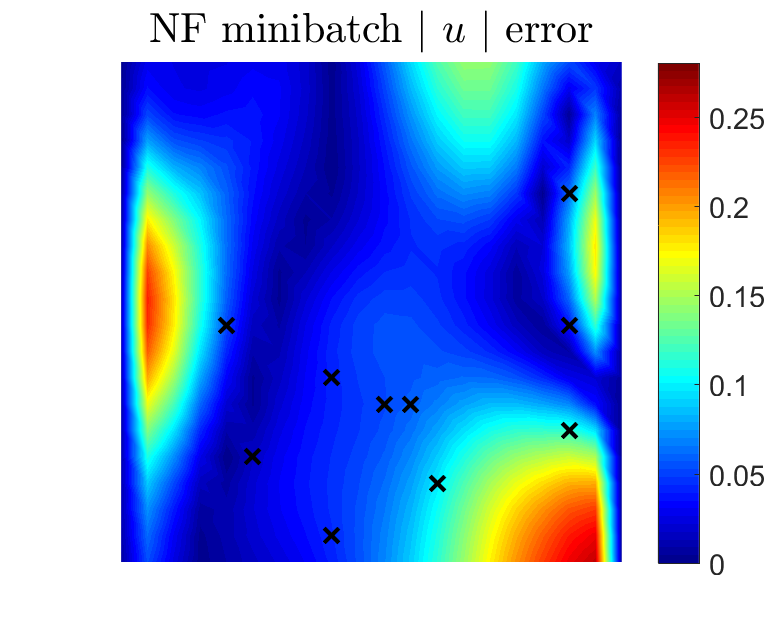}
\includegraphics[width=0.3\textwidth]{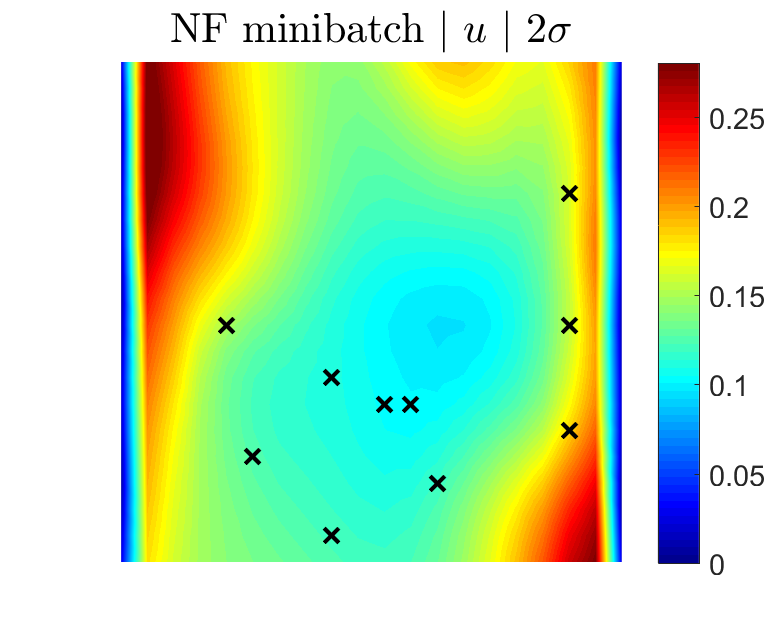}
}
\caption{\label{fig:darcy_nf_mb}
100D Darcy problem: Predictions from NF with minibatch training, i.e., $M_f = 10$, and $M_u = 5$. (a) $\bar{\lambda}$. (b) $u$.  Circle: testing data for $\bar{\lambda}$, Cross: testing data for $u$.
}
\end{figure}

To begin with, we assume that we have  9,900 different paired data  ($\bar{\lambda}$, $u$) as the historical data. For each sample of $\tilde{\lambda}/u$, we utilize $20 \times 20$ uniform grid to resolve it.  We then train a DeepONet to learn the mapping from $\bar{\lambda}(x, y)$ to $u$ using the paired data ($\bar{\lambda}$, $u$). We further train a GAN to learn the functional prior for $\bar{\lambda}(x, y)$ based on the historical data on $\bar{\lambda}(x, y)$. We note that the historical data employed here are the same as in \cite{psaros2023uncertainty,zou2022neuraluq}. 

\begin{figure}
\centering
\subfigure[]{
\includegraphics[width=0.3\textwidth]{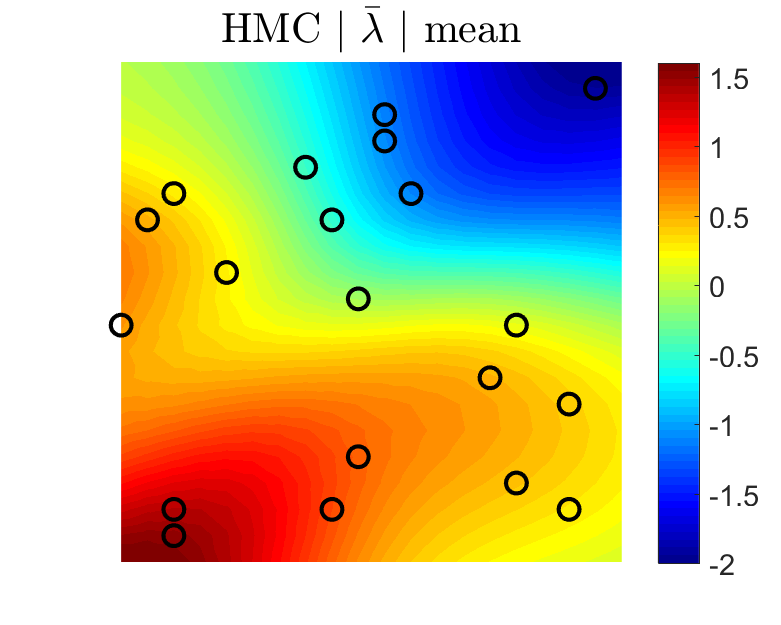}
\includegraphics[width=0.3\textwidth]{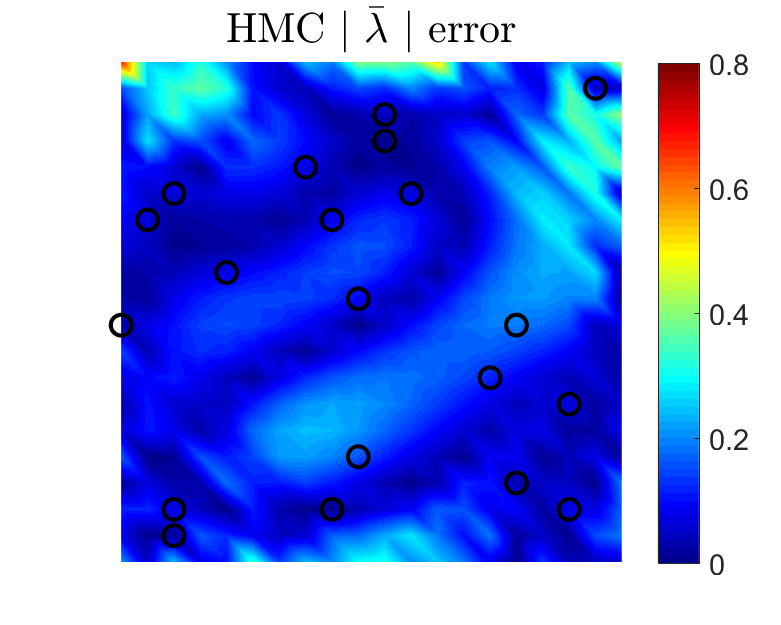}
\includegraphics[width=0.3\textwidth]{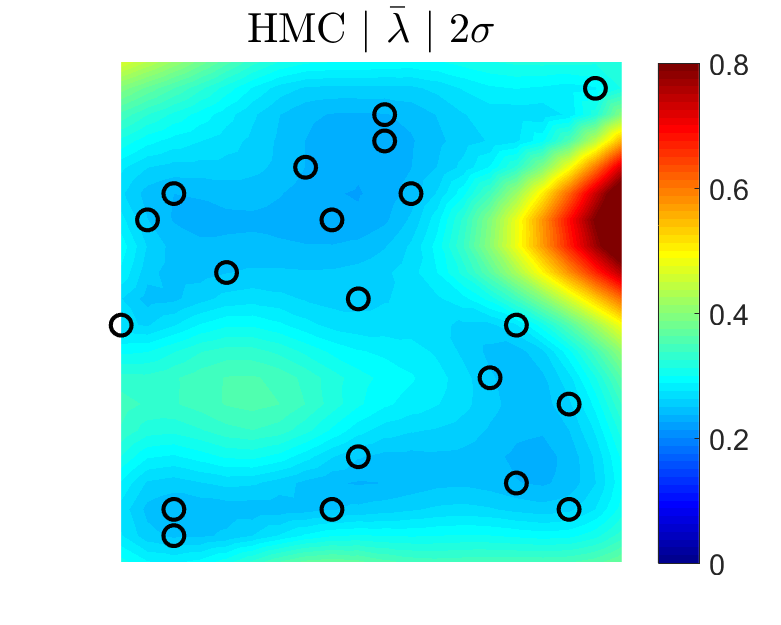}
}
\subfigure[]{
\includegraphics[width=0.3\textwidth]{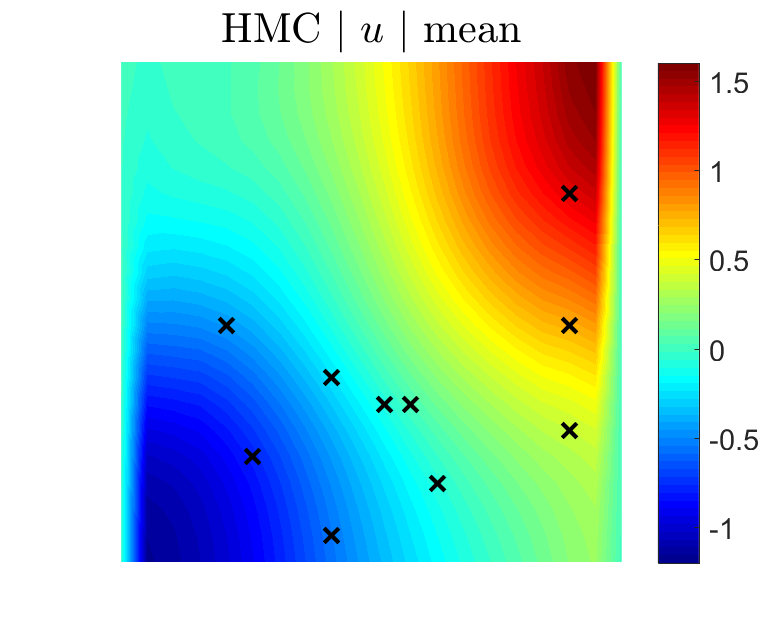}
\includegraphics[width=0.3\textwidth]{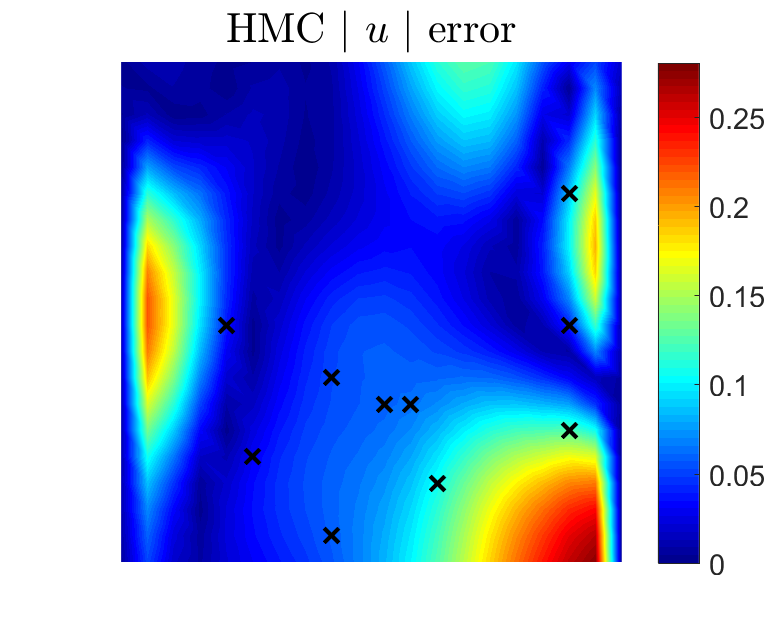}
\includegraphics[width=0.3\textwidth]{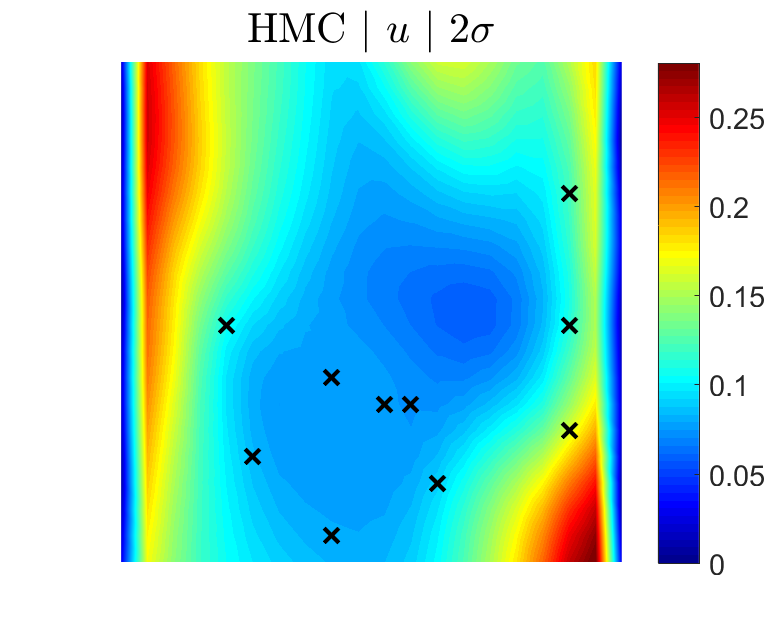}
}
\caption{\label{fig:darcy_hmc}
100D Darcy problem: Predictions from HMC. (a) $\bar{\lambda}$. (b) $u$.  Circle: testing data for $\bar{\lambda}$, Cross: testing data for $u$.
}
\end{figure}

With the pretrained DeepONet and neural functional prior for $\bar{\lambda}$, we assume that we have 20 and 10 noisy measurements for $\lambda$ as well as $u$ for an unseen task, which are displayed in Figs. \ref{fig:darcy_nf_fb}-\ref{fig:darcy_hmc}. The objective is to reconstruct complete $\bar{\lambda}$ and $u$. Similarly, the measurements in the new task are the same as employed in \cite{psaros2023uncertainty,zou2022neuraluq}. We present the predictions for $\bar{\lambda}$ and $u$ from NF with full- and mini-batch training in Figs.~\ref{fig:darcy_nf_fb} and \ref{fig:darcy_nf_mb}, respectively.
It is  observed that the (1) the results in these two cases show little discrepancy, and (2) the computational errors between the predicted means and the reference solutions for $\bar{\lambda}$ and $u$ are mostly bounded by the predicted uncertainties in both cases. We also present the results from HMC for reference in Fig. \ref{fig:darcy_hmc}. As shown, the NF with full- and mini-batch training are able achieve simialr accuracy as HMC in terms of the predicted mean and uncertainties. 


\section{Summary}
\label{sec:summary}
In the present study, we utilize the generative adversarial networks (GANs) to learn the functional prior from historical data and available physics. In addition, two different scenarios for encoding the physics in GANs are considered, i.e., (1) the differential operators for defining the problems are known, we then encode the physical laws via automatic differentiation similar as in PINNs; and (2) the differential operators for defining the problems are unknown, we thus employ the DeepONet to learn the operators given data. We refer to the pre-trained GANs  as neural functional prior.  Further, we propose to employ the normalizing flows to compute the posterior in the latent space of neural functional prior in the context of variational inference as the alternative to the ``gold rule'' HMC. Specially, the NF is a unified framework for both full- and mini-batch training.  We begin with a one-dimensional example to show that NF with full- and mini-batch training are able to achieve similar accuracy comparing to the ``gold rule'' HMC in posterior estimation. 
We further tested 1D and 2D differential equation problems using automatic differentiation and DeepONets to encode physics, respectively. In these two problems, we show that NF can provide accurate predicted means and reasonable uncertainty bounds, with relatively small number of sensors, which can be attributed to the informative functional priors that reflect our knowledge from historical data. Also, NF with full- and mini-batch training is capable of providing similar results as compared to HMC. The minibatch training in NF makes it a promising tool for quantifying uncertainties in high-dimensional parametric PDEs with big data.

\section*{Acknowledgements}
X. Meng would like to acknowledge the support of the National Natural Science Foundation of China (No. 12201229), and the CCF-Baidu Open Fund. X. Meng thanks Dr. Liu Yang for the helpful discussion.


\appendix
\section{Brief overview of DeepoNets}
\label{sec:deeponet}
\begin{figure}
\centering
\includegraphics[width=0.45\textwidth]{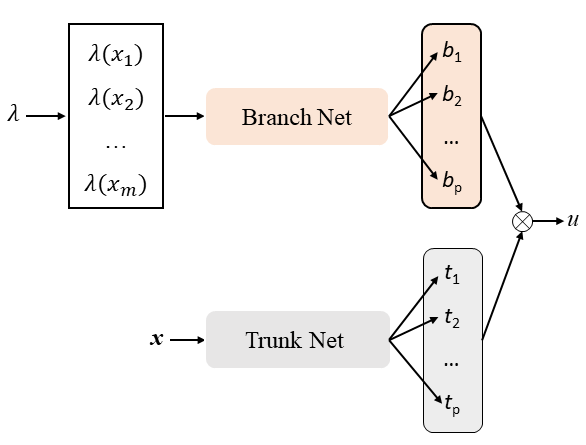}
\caption{\label{fig:deeponet}
Schematic of DeepONet. $\lambda$ is the input of  Branch Net,  $x_1, x_2, .., x_m$ are the discrete points to resolve the input function, $b_1, b_2, ..., b_p$ and $t_1, t_2, ..., t_p$ are the outputs for the Branch Net and Trunk Net, respectively, ${u}$ is the output target function. If the operator takes multiple functions as input (e.g., source term $f$, boundary condition $b$, etc.), then we just need to concatenate the multiple vectors that resolve these functions as the input of the Branch Net. Adapted from \cite{meng2022learning}.
}
\end{figure}

 As reported in \cite{lu2021learning},  DeepONets  can be used as a universal approximator to any continuous nonlinear operator. We present a schematic of DeepONets in Fig. \ref{fig:deeponet}.  As shown, the DeepONet has two sub-networks, i.e., the Branch Net (BN) and the Trunk Net (TN). The input for BN is a function, which is represented by a set of discrete function values at certain locations, i.e., $x_1, x_2, ..., x_m$, and the output of BN is a vector $[b_1, b_2, ..., b_p]$. Further,  TN takes $\bm{x}$ as input and outputs a vector $[t_1, t_2, ..., t_p]$. The output of the DeepONet is the inner product of these two vectors as $\bm{u} = \sum^p_{j=1} b_j t_j$. 
 The DeepONet is essentially a mapping between two function spaces, e.g., from the problem parameters $\lambda$ to the solution of a PDE, i.e., $u$.

 Similar as in \cite{meng2022learning}, the uniform grids are used to  discretize the input functions in the test cases of the present study. Also, we do not employ any constraint on the input $\bm{x}$ in the TN, we can thus evaluate the output ${u}$ at any location. Both the BN and TN are trained simultaneously by minimizing  the mean squared error (MSE) between the given and predicted ${u}$ from DNNs using the Adam optimizer. More details on DeepONets can be found in \cite{lu2021learning,lu2022comprehensive}.

\section{Details for numerical computations}
\label{sec:nndetails}
In this section, we present the details for  data generation, architectures of employed NNs, as well as the training strategy.

For the training data, the historical data as well as the testing data used in Secs. \ref{sec:poisson} and \ref{sec:darcy} are the same as those in \cite{meng2022learning} and \cite{psaros2023uncertainty,zou2022neuraluq}, respectively, which are available on \href{https://github.com/Crunch-UQ4MI}{\textit{github.com/Crunch-UQ4MI}}.

The architectures of the GANs for neural functional priors and NFs employed in the present work are displayed in Tables \ref{table:gan_arch} and \ref{table:nf_arch}, respectively. Note that the architectures of GANs are the same as in \cite{meng2022learning}. We employ the Adam with a learning rate $10^{-4}$ for both the training of neural functional priors as well as NFs.  Furthermore, for the HMC, we employ the No-U-Turn \cite{hoffman2014no}, which can adaptively set path lengths in HMC in this study. In all cases, the initial step size is set as 1, the target acceptance rate is 0.6, and the number of burnin steps is 2,000. Particularly, the No-U-Turn is implemented using the Tensorflow Probability package \cite{lao2020tfp}. 

\begin{table}[H]
\centering
{\footnotesize
\begin{tabular}{c|cc|cc|c|c}
\hline \hline
  & \multicolumn{2}{c|}{G ($\tilde{g}/\tilde{h}^*$)}  & \multicolumn{2}{c|}{D} & \multirow{2}{*}{$d_{\bm{\xi}}$} & \multirow{2}{*}{Training steps}  \\
  & width $\times$ depth  & Activation & width $\times$ depth & Activation &\\
  \hline
  
  {Sec. \ref{sec:func}} & $64 \times 2 / 64 \times 2$   & tanh/tanh & $128 \times 3$ & Leaky ReLu &   10 &500,000 \\
  \hline
  {Sec. \ref{sec:poisson}} & $64 \times 2 / 64 \times 2$   & tanh/tanh & $128 \times 3$ & Leaky ReLu &  40 & 500,000 \\
  \hline
{Sec. \ref{sec:darcy}} & $64 \times 2 / 64 \times 2$   & tanh/tanh & $128 \times 3$ & Leaky ReLu &  100 & 500,000 \\
  \hline \hline
\end{tabular}
}
\caption{
Architecture and training steps of GANs in each case. The width and depth are for the hidden layers. $d_G$ is the dimension of $\bm{\xi}$ as well as the output dimension of $\tilde{g}$ and $\tilde{h}^*$ in $G$.
}
\label{table:gan_arch}
\end{table}

\begin{table}[H]
\centering
{\footnotesize
\begin{tabular}{c|c|c|c|c}
\hline \hline
  &bolcks & width $\times$ depth  & Activation &  Training steps\\
  \hline
  
  {Secs. \ref{sec:func}-\ref{sec:darcy}} & 4 & $256 \times 2 $   & tanh &  200,000 \\
  \hline \hline
\end{tabular}
}
\caption{
Architecture and training steps of NF in each case. The width and depth are for the hidden layers in each block of NF.
}
\label{table:nf_arch}
\end{table}





\bibliographystyle{elsarticle-num}
\biboptions{sort&compress}
\bibliography{refs.bib}







\end{document}

%% file: math_commands.tex
\usepackage{amsmath,amsfonts,bm}

\def\eqref#1{equation~(\ref{#1})}
\def\Eqref#1{Equation~(\ref{#1})}

\def\1{\bm{1}}

\def\vx{{\bm{x}}}

\DeclareMathAlphabet{\mathsfit}{\encodingdefault}{\sfdefault}{m}{sl}
\SetMathAlphabet{\mathsfit}{bold}{\encodingdefault}{\sfdefault}{bx}{n}

